\documentclass{easychair}
\usepackage{amsfonts, amsthm, amssymb, amsmath, stackengine, scalerel}
\usepackage{mathrsfs,array,tikz-cd}
\usepackage{hyperref}
\usepackage{eucal,fullpage,times,color,enumerate,accents, comment}
\usepackage[all]{xy}
\usepackage{xr}
\usepackage{url}
\usepackage{enumitem}
\usepackage[new]{old-arrows}
\usepackage{extpfeil}
\usepackage{turnstile} 
\usepackage{verbatim}  
\usepackage{graphics,float} 
\usetikzlibrary{graphs,decorations.pathmorphing,decorations.markings}
\input xy
\xyoption{all}
\usepackage{adjustbox}

\usepackage{algorithm,caption}  
\DeclareCaptionFont{lightblue}{\color{teal}}
\newextarrow{\xbigtoto}{{20}{20}{20}{20}}
{\bigRelbar\bigRelbar{\bigtwoarrowsleft\rightarrow\rightarrow}}

\newcommand{\calU}{{\mathcal{U}}}
\newcommand{\calV}{{\mathcal{V}}}
\newcommand{\calA}{{\mathcal{A}}}
\newcommand{\calE}{\mathcal{E}}
\newcommand{\calB}{{\mathcal{B}}}

\newcommand{\calD}{{\mathcal{D}}}

\newcommand{\calF}{\mathcal{F}}

\newcommand{\calI}{\mathcal{I}}

\newcommand{\calH}{\mathcal{H}}

\newcommand{\calP}{\mathcal{P}}

\newcommand{\Mod}{\mathrm{Mod}}

\newcommand{\R}{\mathbb{R}}

\newcommand{\id}{\mathrm{id}}

\newcommand{\thT}{\mathbb{T}}

\newcommand{\E}{\mathcal{E}}

\DeclareFontFamily{U}{wncy}{}
\DeclareFontShape{U}{wncy}{m}{n}{<->wncyr10}{}
\DeclareSymbolFont{mcy}{U}{wncy}{m}{n}
\DeclareMathSymbol{\Sh}{\mathord}{mcy}{"58}

\DeclareSymbolFont{matha}{OML}{txmi}{m}{it}
\DeclareMathSymbol{\varv}{\mathord}{matha}{118}

\newcommand{\Eff}{\mathrm{Eff}}

\def\forkindep{\mathrel{\raise0.2ex\hbox{\ooalign{\hidewidth$\vert$\hidewidth\cr\raise-0.9ex\hbox{$\smile$}}}}}

\newcommand{\Fin}{\mathrm{Fin}}

\newcommand{\Pw}{\calP(\w)}

\newcommand{\dn}{{\lnot\lnot}}

\newcommand{\lrk}{\leq_{\mathrm{RK}}}
\newcommand{\erk}{\equiv_{\mathrm{RK}}}

\newcommand{\glt}{\leq_{\mathrm{LT}}^{\mathsf{o}}}
\newcommand{\elt}{\equiv_{\mathrm{LT}}^{\mathsf{o}}}
\newcommand{\lT}{\leq_{\mathrm{Tuk}}}
\newcommand{\eT}{\equiv_{\mathrm{Tuk}}}

\newcommand{\w}{\omega}

\newcommand{\clt}{\leq_{\mathrm{LT}}}
\newcommand{\eclt}{\equiv_{\mathrm{LT}}}
\newcommand{\sclt}{<_{\mathrm{LT}}}

\newcommand{\bmid}{\mathrel{\big|}} 

\newcommand{\Dt}{\mathcal{D}_{\rm T}}

\makeatletter
\newcommand*{\relrelbarsep}{.386ex}
\newcommand*{\relrelbar}{%
	\mathrel{%
		\mathpalette\@relrelbar\relrelbarsep
	}%
}
\newcommand*{\@relrelbar}[2]{%
	\raise#2\hbox to 0pt{$\m@th#1\relbar$\hss}%
	\lower#2\hbox{$\m@th#1\relbar$}%
}
\providecommand*{\rightrightarrowsfill@}{%
	\arrowfill@\relrelbar\relrelbar\rightrightarrows
}
\providecommand*{\leftleftarrowsfill@}{%
	\arrowfill@\leftleftarrows\relrelbar\relrelbar
}
\providecommand*{\xrightrightarrows}[2][]{%
	\ext@arrow 0359\rightrightarrowsfill@{#1}{#2}%
}
\providecommand*{\xleftleftarrows}[2][]{%
	\ext@arrow 3095\leftleftarrowsfill@{#1}{#2}%
}
\makeatother

\makeatletter
\newcommand*\Lslash
{%
	\text{\rlap{\kern.03em\@xxxii}}L%
}
\makeatother

\newcommand{\cosimp}[3]{\xymatrix@1{#1 \ar@<.4ex>[r] \ar@<-.4ex>[r] & {\ }#2 \ar@<0.8ex>[r] \ar[r] \ar@<-.8ex>[r] & {\ } #3 \ar@<1.2ex>[r] \ar@<.4ex>[r] \ar@<-.4ex>[r] \ar@<-1.2ex>[r] & \cdots }}
\newsavebox{\pullback}
\sbox\pullback{%
	\begin{tikzpicture}%
	\draw (0,0) -- (1ex,0ex);%
	\draw (1ex,0ex) -- (1ex,1ex);%
	\end{tikzpicture}}

\makeatletter
\newcommand*\dotp{\mathpalette\dotp@{.5}}
\newcommand*\dotp@[2]{\mathbin{\vcenter{\hbox{\scalebox{#2}{$\m@th#1\bullet$}}}}}
\makeatother

\newcommand{\equalizer}[2]{\xymatrix@1{#1 \ar@<.4ex>[r] \ar@<-0.4ex>[r] & {\ } #2}}

\newcommand{\adjunction}[4]{\xymatrix@1{#1{\ } \ar@<-0.3ex>[r]_{ {\scriptstyle #2}} & {\ } #3 \ar@<-0.3ex>[l]_{ {\scriptstyle #4}}}}

\usepackage{epigraph}

\usepackage{xcolor}
\usepackage{csquotes}
\usepackage{lipsum}

\definecolor{quotemark}{gray}{0.7}
\makeatletter
\newlength\origparskip

\newcommand{\fquote}{%
	\@ifnextchar[{\fquote@i}{\fquote@i[]}
}

\def\fquote@i[#1]{%
	\@ifnextchar[{\fquote@ii{#1}}{\fquote@ii{#1}[]}
}%

\def\fquote@ii#1[#2]{%
	\def\pqm@tempa{#1}%
	\def\pqm@tempb{#2}%
	\noindent
	\list
	{}
	{\setlength{\leftmargin}{0.3\textwidth}%
		\setlength{\rightmargin}{0.1\textwidth}%
		\setlength{\origparskip}{\parskip}}%
	\item[]%
	\begin{picture}(0,0)%
	\put(-15,-8){\makebox(0,0){\scalebox{4}{%
				\textcolor{quotemark}{\textquotedblright}}}}%
	\end{picture}%
	\begingroup
	\itshape
	\ignorespaces}%

\def\endfquote{%
	\endgroup
	\par
	\raggedleft
	\ifx\pqm@tempa\empty
	\else
	{\bfseries --- \pqm@tempa\par}%
	\setlength{\parskip}{\origparskip}%
	\ifx\pqm@tempb\empty
	\else
	(\pqm@tempb)%
	\fi
	\fi
	\par
	\endlist}
\makeatother

\begin{document}
	\bibliographystyle{alpha}
	\newtheorem{theorem}{Theorem}[section]
	\newtheorem*{theorem*}{Theorem}
	\newtheorem*{condition*}{Condition}
	\newtheorem*{definition*}{Definition}
	\newtheorem*{corollary*}{Corollary}
	\newtheorem{proposition}[theorem]{Proposition}
	\newtheorem{lemma}[theorem]{Lemma}
	\newtheorem{corollary}[theorem]{Corollary}
	\newtheorem{claim}[theorem]{Claim}
	\newtheorem{conclusion}[theorem]{Conclusion}
	\newtheorem{hypothesis}[theorem]{Hypothesis}
	\newtheorem{conjecture}[theorem]{Conjecture}
	\newtheorem{setup}[theorem]{Setup}
	\newtheorem{sumthm}[theorem]{Summary Theorem}

	\newtheorem{maintheorem}{Theorem}
	\renewcommand*{\themaintheorem}{\Alph{maintheorem}}
	\newtheorem{mainprop}[maintheorem]{Proposition}
	
	\theoremstyle{definition}
	\newtheorem{definition}[theorem]{Definition}
	\newtheorem{question}[theorem]{Question}
	\newtheorem{action}[theorem]{Action Item}
	\newtheorem{answer}[theorem]{Answer}
	\newtheorem{goal}[theorem]{Goal}
	\newtheorem{exercise}[theorem]{Exercise}
	\newtheorem{remark}[theorem]{Remark}
	\newtheorem{observation}[theorem]{Observation}
	\newtheorem{discussion}[theorem]{Discussion}
	\newtheorem{guess}[theorem]{Guess}
	\newtheorem{example}[theorem]{Example}
	\newtheorem{condition}[theorem]{Condition}
	\newtheorem{warning}[theorem]{Warning}
	\newtheorem{notation}[theorem]{Notation}
	\newtheorem{construction}[theorem]{Construction}
	
	\newtheorem{problem}{Problem}
	\newtheorem{fact}[theorem]{Fact}
	\newtheorem{thesis}[theorem]{Thesis}
	\newtheorem{convention}[theorem]{Convention}
	\newtheorem{summary}[theorem]{Summary}


	%
	\title{What can Topology tell us about Logical Complexity?}

	\author{Takayuki Kihara\thanks{TK was partially supported by JSPS KAKENHI Grant Numbers 22K03401 and 23K28036.} \,  and Ming Ng\thanks{MN was partially supported by a JSPS PostDoc Fellowship (Short-Term).}}
	
	%
	%
	\institute{
		Nagoya University,
		Japan \\
		\email{kihara@i.nagoya-u.ac.jp}\\
		\email{ngming.math@gmail.com}
	}
	
	
	\authorrunning{Kihara and Ng}
	
	
	\titlerunning{What can Topology tell us \dots }

	\maketitle
	
	\begin{abstract} In the 1980s, category theorists introduced the Lawvere-Tierney $(\leq_{\mathrm{LT}})$ order in the Effective Topos, known to effectively embed the Turing degrees. Understanding its structure is a longstanding open problem in the area.  In particular, there was an informal sense that the $\leq_{\mathrm{LT}}$-order reflects certain shifts in combinatorial complexity, but a precise characterisation remained elusive for some time.
		
		Recent work by the authors has substantially clarified the picture. In \cite{KiNg26}, the authors introduced a game-theoretic (``gamified'') version of the Kat\v{e}tov order on filters over $\omega$ -- essentially, this is the usual Kat\v{e}tov order now closed under well-founded iterations of Fubini powers. The first major theorem of the paper was to show that a computable variant of the gamified Kat\v{e}tov order is isomorphic to the original $\leq_{\mathrm{LT}}$-order. This was a surprising discovery, and opens up many challenging questions regarding the interplay between combinatorial and computable complexity, which informed the rest of the paper's investigations.
		
		This note gives an informal survey of some of these interactions explored in \cite{KiNg26}, and announces some forthcoming results. The guiding perspective is that different notions of complexity arising in different areas of logic can be seen to be controlled by the same mechanism -- once placed in the right topological framework.
	\end{abstract}

	\section{What is a generalised space?}
	
	Logicians are familiar with the key notion of a filter, either as a generalised point or as an abstract notion of largeness. Topos theorists push this idea much further by asking: what is a {\em generalised space}?
	
	Our main protagonist is the {\em Lawvere-Tierney topology} [hereafter: LT topology]. Informally, an LT topology gives an abstract notion of localness. A space $X$ might be said to satisfy some property $P$ locally if $X$ can be covered by a family of opens  $\{U_i\}_{i\in I}$, each satisfying the property $P$. Abstracting this, one might require:
	\begin{enumerate}
		\item If ``$P$ implies $Q$'', then ``locally-$P$ implies locally-$Q$''
		\item $P$ implies locally-$P$
		\item Locally-locally-$P$ is equivalent to locally-$P$.
	\end{enumerate}

Such notions can be defined in any sufficiently well-behaved category (i.e. a topos), and admit several equivalent formulations -- see e.g. \cite[\S V.1]{MM} or \cite[\S 16]{HylandEffective}. Our focus is the {\em Effective Topos} $\Eff$, a category-theoretic framework for computability, where LT topologies take the following form.

	
	\begin{definition}\label{def:LT-topology} 
		An {\em LT topology} in $\Eff$ is an endomorphism 
		$$j\colon \Pw\to \Pw$$  
		subject to the following conditions:
		\begin{enumerate}
			\item $\forall p,q\in\Pw. \left( p\to q\right )\to \left(j(p)\to j(q)\right)$
			\item  $\forall p\in\Pw. \left( p\to j(p)\right)$
			\item  $\forall p\in\Pw.\,\, j\circ j(p)=j(p)$
		\end{enumerate}	
		Here, $\Pw$ denotes the powerset of the natural numbers, and ``$\to$'' denotes the standard realisability interpretation of implication. That is, fixing an enumeration of partial computable functions $\{\varphi_{e}\}_{e\in\w}$,
		$$p\to q:=\left\{e\mid  \forall a\in p\,,\, \varphi_{e}(a) \,\text{is defined and}\, \varphi_e(a)\in q\right\}.$$
For instance, Condition (2) asserts the existence of some $e\in\w$ whereby
$e\in (p\to j(p))$ for all $p\in\Pw$.
	\end{definition} 
	
	Extending this perspective, we define the following:
	
	\begin{definition} In any topos $\calE$, the set of LT topologies in $\calE$ carries a natural pre-order called the {\em Lawvere-Tierney order} [hereafter: $\clt$-order]. When $\calE$ is the Effective Topos, this pre-order is defined as: 
		$$j\clt k : \iff \forall p.j(p)\to k(p) \, \text{is valid}\,.$$
	\end{definition}
	
	A key structural fact is that the Turing Degrees embed effectively into this $\clt$-order \cite[\S 17]{HylandEffective}, which indicates that the $\clt$-order gives a measure of computable complexity. Returning to the informal picture we gave before, one may ask: if LT topologies correspond to generalised spaces, what then are its (generalised) points? Moreover, how do they interact with the $\clt$-order?
	
	The situation is well-understood whenever the topos $\calE$ is of topological origin, such as $\calE=\textbf{Sh}(\R)$, the category of sheaves on $\R$. In which case, the global points of $\textbf{Sh}(\R)$ correspond to completely prime filters in the lattice of open sets in $\R$, and the LT topologies in $\textbf{Sh}(\R)$ correspond to subspaces of $\R$ -- see \cite[\S 2 and \S 7.2]{Vi07}. By contrast, in more abstract settings such as $\E=\Eff$, the situation is much less well understood. In fact, progress on understanding the $\clt$-order in $\Eff$ has been relatively slow. Early work by Hyland and Pitts \cite{HylandEffective,PittsPhD} established the existence of a minimum and maximum class on the non-trivial LT topologies. Attention therefore turned to investigating what lay in between the two extreme cases, a longstanding open problem in the area:
	$$\id \sclt \,\,\cdots\cdots \,? \,\cdots\cdots \,\, \sclt \dn\,\,.$$

	\section{Towards the Gamified Kat\v{e}tov order}

	A key turning point emerged in Lee-van Oosten's striking discovery that all LT topologies in $\Eff$ are constructed from certain basic building blocks \cite{Lee,LvO13}. The following summary theorem makes this statement precise, and includes a couple other highlights. 
	\begin{sumthm}[{{\cite{LvO13}}}]\label{sumthm-1} Any family of subsets $\calA\subseteq \Pw$ defines an LT topology, which we denote as $$j_\calA\colon \Pw\to \Pw.$$
		An LT topology of this form is called a {\em basic topology}. In particular:
		\begin{enumerate}[label=(\roman*)]
			\item Every LT topology in $\Eff$ is a recursive join of basic topologies.
			\item Denote $\bigcap \calA$ to be the intersection of all $A\in\calA$. \underline{Then}: 
			$$j_\calA=\id \iff \bigcap\calA\neq \emptyset\quad .$$
			\item Say $\calA$ has the {\em $n$-intersection property} if for any $n$ elements $A_1,\dots,A_n\in\calA$, their intersection is non-empty: $$\bigcap^{n}_{i=1} A_i\neq \emptyset.$$ 
			Suppose $\calA\subseteq\Pw$ has the $n$-intersection property whereas $\calB\subseteq\Pw$ {\em fails} this property. 
			\underline{Then}:
			$$j_\calB \not\clt j_{\calA}.$$
		\end{enumerate}
	\end{sumthm}
	
	The crux move of our paper \cite{KiNg26} was to regard the $\clt$-order as inducing a pre-order on subset families $\calA\subseteq \Pw$. In which case, item (ii) of the Summary Theorem has the following elegant translation: $\calA$ belongs to the $\eclt$-minimum class iff $\bigcap\calA\neq\emptyset$. This restatement will be suggestive to the set theorist because of the following well-known fact: 
	
	\begin{fact}[see e.g. {{\cite[Exercise 1.6.3]{Gold22}}}] The {\em Rudin-Keisler order} ($\lrk$) defines a preorder on ultrafilters over $\w$. In particular, there exists a minimum class, where $\calF\subseteq\Pw $ belongs to the $\erk$-minimum class iff $\bigcap\calF\neq\emptyset$.
	\end{fact}
	
	As it turns out, the Rudin-Keisler order is well-defined not just on ultrafilters, but also on upper sets over $\w$ (i.e. subset families closed upwards under $\subseteq$). This sets up the first major theorem of \cite{KiNg26}: 
	
	\begin{theorem}[{{\cite[Theorem A]{KiNg26}}}]\label{thm:mainthm}  The Gamified Kat\v{e}tov order, written suggestively as
		$$\calU\glt\calV,\qquad\qquad\qquad\calU,\calV\subseteq \; \text{upper sets}\,,$$
		is a preorder on upper sets over $\w$. In particular:
		\begin{enumerate}[label=(\roman*)]
			\item The \emph{Gamified Kat\v{e}tov order} is equivalent to the {\em Kat\v{e}tov order} closed under {\em well-founded iterations of Fubini powers}. In particular, it is strictly coarser than {\em Rudin-Keisler} and {\em Kat\v{e}tov}.
			\item The {\em Gamified Kat\v{e}tov order} admits an explicit game-theoretic description, justifying its name.
			\item The {\em computable Gamified Kat\v{e}tov order} is equivalent to the {\em $\clt$-order} on upper sets over $\w$.
		\end{enumerate}
	\end{theorem}
	
	\begin{remark}\label{rem:LT} In the same paper, we later extend the Gamified Kat\v{e}tov order from upper sets to {\em upper sequences} (i.e. a countable sequence of upper sets), and show that the computable version of this extended Gamified Kat\v{e}tov order is isomorphic to the original $\clt$-order \cite[Theorem 6.5]{KiNg26}. This move from upper sets to upper sequences mirrors the move from basic topologies to recursive join of basic topologies, as in item (i) of the Summary Theorem.
	\end{remark}

	Theorem~\ref{thm:mainthm} highlights several design features that distinguish the Gamified Kat\v{e}tov order from other preorders in the set theory literature. Taken by themselves, these choices may seem puzzling to the reader, but in fact they reflect a deep structural phenomenon, emerging in three different guises: 
	
	\begin{itemize} 
		\item[$\diamond$] as a natural variant of the Kat\v{e}tov order (set theory);
		\item[$\diamond$] as a concrete representation of abstract topologies, namely LT topologies (category theory);
		\item[$\diamond$] as a new hierarchy of computability notions, namely ``computability by majority''\footnote{This involves running multiple computations in parallel, and adopting the result decreed by the ``majority'' of these computations; for a discussion on what majority means and further details, see \cite[\S 7.2]{KiNg26}.} (computability theory).
	\end{itemize}
	Seen in this light, the topos theory brings to the surface a series of deep structural connections we would otherwise be hard-pressed to identify. With Theorem~\ref{thm:mainthm} in place, we now have the vocabulary to ask many interesting new questions, laying the groundwork for a wide-ranging research programme. 
	
	\section{``Coarse and Subtle''}
	
	One may view the $\glt$-order as simply isolating the combinatorial mechanism underlying the original $\clt$-order, but it is also interesting in its own right. Theorem~\ref{thm:mainthm} states that the $\glt$-order is strictly coarser than the Rudin-Keisler order, but the nature of its coarseness is subtle. Preliminary evidence indicates that the Gamified Kat\v{e}tov order detects a kind of complexity quite different from that detected by the classical set-theoretic tools. For instance: 
	
	\begin{theorem}[{{\cite[Theorem B]{KiNg26}}}]\label{thm:Tukey} The Gamified Kat\v{e}tov and Tukey order are incomparable on filters over $\w$ (in ZFC).
	\end{theorem}
	
	It is well-known that the Tukey order ($\lT$) is also strictly coarser than the Rudin-Keisler order. One may therefore wonder if these two orders are in fact responding to the same kind of structural properties of filters. Theorem~\ref{thm:Tukey} shows this is not the case: the Gamified Kat\v{e}tov and Tukey orders are in fact fundamentally misaligned. In particular, the $\glt$-order does \emph{not} measure the cofinality types of filters. 
	
	Moreover, the fact that Theorem~\ref{thm:Tukey} holds in ZFC is also noteworthy. A recent  breakthrough result of \cite{CaZa24} shows that it is consistent with ZFC that all non-principal ultrafilters on $\w$ are Tukey-equivalent. By contrast, part of the proof of Theorem~\ref{thm:Tukey} shows that each $\elt$-class contains at most $2^{\aleph_0}$ ultrafilters in ZFC. This rigidity was then leveraged to deduce the existence of ultrafilters $\calU, \calV$ such that $\calU \eT \calV$ but $\calU \not\elt \calV$.
	
	\medskip
	If the $\glt$-order does not measure the cofinality types of filters, what kind of combinatorial complexity does it actually measure? One way to investigate this is to compare its structure with the classical Kat\v{e}tov order, and examine what the $\glt$-order regards as significant vs. insiginficant differences between filters (dually, ideals) on $\w$. 
	\smallskip 
	
	As a warm-up, recall that an AD family (= {\em Almost Disjoint}) is a subset family $\calA\subseteq\Pw$ such that all member sets are infinite and have finite pairwise intersection. A MAD family (= {\em Maximal Almost Disjoint}) is an AD family $\calA$ that is maximal with respect to this property. Informally, MAD families are very flat -- their member sets are spread out across $\w$ such that they only overlap in a small area. 
	
	Each MAD family generates an ideal on $\w$. From one point of view, these MAD families can be seen as mapping out a wide range of complexity. Within the classical Kat\v{e}tov order, there exists at least $\mathfrak{c}^+$-many equivalence classes of MAD families as well as  $\mathfrak{c}$-many pairwise incomparable ones \cite[\S2]{HruGF03}. By contrast, in the gamified setting, all MAD families collapse to a single $\elt$-equivalence class \cite[Conclusion 5.3]{KiNg26}. On the other hand, despite its apparent coarseness, there still remains a rich internal structure within $\glt$. A first indication of this is \cite[Theorem C]{KiNg26}, which establishes the existence of an infinite strictly ascending chain of classes within the Gamified Kat\v{e}tov order. In a forthcoming paper, we significantly extend this picture:
	
	\begin{theorem}[{{\cite{KiNg-nonlinear}}}]\label{thm:Pwfin} The Gamified Kat\v{e}tov order admits an embedding of $\Pw/\Fin$. In particular, it contains an increasing chain of length $\mathfrak{b}$, as well as an antichain of size $\mathfrak{c}$.
	\end{theorem}
	
	Let us step back for a moment. Given what we now know about the Gamified Kat\v{e}tov order, how does this illuminate the structure of the original $\clt$-order? What can we say about the interactions between the combinatorial and computable within this framework? 
	
	These question guide the work of Sections 5-6 in \cite{KiNg26}. When the terrain is still relatively uncharted, it is difficult to make definitive statements on where the path may take us. Still, a couple of interesting takeaways have emerged that are worth highlighting to the general logical audience. 
	
	\smallskip 
	First, a key theme in our recent investigations is that if we wish to separate equivalence classes within the $\clt$-order, it is helpful to divide the analysis into two stages: start by performing the separation in the (non-computable) setting of the Gamified Kat\v{e}tov order, before re-introducing the computability constraints to obtain the final separation in $\clt$. In particular, since $\glt$ is coarser than $\clt$ (a consequence of Theorem~\ref{thm:mainthm}), this means that $\calH\not\glt\calI$ implies $\calH\not\clt\calI$. This insight was applied in \cite[Corollary 5.7]{KiNg26}, which showed that given the infinite increasing chain in $\glt$ established by \cite[Theorem C]{KiNg26}, one obtains the same infinite chain in $\clt$ essentially for free. 
	\smallskip 
	
	Second, it is worth remarking that serious connections between set theory and computability theory are, of course, not new.  What is distinctive, however, in the present setting is that the $\clt$-order provides a uniform framework in which both forms of complexity can be compared directly. For instance, applying Summary Theorem~\ref{sumthm-1}, every filter $\calF\subseteq\Pw$ determines an initial segment of the Turing degrees:
	$$
	\calD_{\rm T}(\calF):=\big \{ \,[f\colon\w\to\w] \bmid f\clt \calF \big\}.
	$$
	This defines an interesting new invariant of filters, which we call the \emph{Turing degree profile} of $\calF$.
	
	It is reasonable to wonder if filters and Turing degrees are largely orthogonal within $\clt$, with relatively few direct points of contact. Encouragingly, the following Cofinality Theorem shows this is very much not the case.
	
	\begin{theorem}[{{Cofinality, \cite[Theorem E]{KiNg26}}}]\label{thm:cofinality} The following is true:
		\begin{enumerate}[label=(\roman*)]
			\item For any $f\colon\omega\to\omega$, there exists a summable ideal $\calI$ such that
			$$f \clt \calI^*,$$
			where $\calI^*$ is the dual filter.
			\item There does not exist a filter $\calF\subseteq\Pw$ such that 
			$$f\clt\calF,\qquad\text{for all Turing oracles}\; f\;.$$
		\end{enumerate}
		
	\end{theorem}

	\noindent In English: every Turing degree $f$ appears in the profile of some filter $\calF$, but no single filter captures them all. Theorem~\ref{thm:cofinality} therefore identifies an unexpected connection between the combinatorial and the computable: one can systematically probe the structure of filters on $\w$ by examining which Turing degrees they bound, and which they necessarily omit.

	As a first step towards developing this perspective, we characterise the Turing degree profiles of all non-principal $\Delta^1_1$-filters, significantly generalising previous work of van Oosten \cite{vO14} and the first author \cite{Kih23}: 
	
	\begin{theorem}[{{\cite[Theorem F]{KiNg26}}}] Let $\calF\subseteq \Pw$ be a non-principal $\Delta^1_1$-filter. Then, for any $f\colon\w\to\w$, 
		$$f\clt \calF \iff f\;\;\text{is hyperarithmetic}\;.$$
	\end{theorem}
	
	\begin{remark} Notice the definition of the Turing degree profile is quite flexible. If $\Dt(\calF)$ turns out to be too coarse for one's purpose, one can always analyse the filters based on, e.g. their Weihrauch degree profiles etc. 
	\end{remark}
	
	\section{Some Final Remarks}
	
	The Effective Topos $\Eff$ was introduced by Hyland in the early 80s \cite{HylandEffective}, whereas the Rudin-Keisler order dates back to the early 70s (if not earlier) \cite{Rud71,BlassThesis}. Certainly the importance of both constructions were independently well-appreciated, and have been extensively studied. Why, then, did it  take forty-odd years after Hyland's seminal paper for the connection to finally be made? 
	
	One natural explanation lies in the historically limited interactions between category theory and set theory, but there is more to be said here. The Gamified Kat\v{e}tov order emerged from re-examining previous results by Lee-van Oosten \cite{LvO13} (who showed the $\clt$-order reflects differences in intersection patterns of subset families) and by the first author \cite{Kih23} (who showed the $\clt$-order admits game-theoretic analysis). Looked at in this light, it seems that certain structural results about the $\clt$-order needed to be in place before one might even begin formulating a result like Theorem~\ref{thm:mainthm}.

	
	Moreover, recognising the possibility of a connection is only part of the story. As noted in the abstract of \cite{KiNg26}: ``The proofs [of the paper] draw on ideas from general topology, descriptive set theory, and computability theory.'' In particular, none of the main arguments in \cite{KiNg26} use any category theory --- despite the category-theoretic origins of the $\clt$-order. In this sense, establishing the connection required a genuine crossing of mathematical boundaries. Perhaps then it is not so surprising that the link remained hidden for so long: there was little {\em a priori} reason to expect that combinatorial set theory would provide an important key to understanding the $\clt$-order.
	
	After we proved Theorem~\ref{thm:mainthm}, many other results soon followed -- surely there is still much more to be said. We refer the interested reader to \cite[\S 7]{KiNg26}, the open problems section of that paper. More broadly, our results point to an emerging new understanding of category theory,  particularly in regards to its relationship with logical complexity. Beyond just providing an abstract framework for mathematics, category theory also implicitly determines which distinctions should be regarded as structurally meaningful, and which should be treated as inessential. In categorical logic, for instance, two first-order theories $\thT_0,\thT_1$ may fail to be bi-interpretable yet still be regarded as equally expressive in the sense that $\Mod(\thT_0)\simeq \Mod(\thT_1)$, i.e. their categories of models are equivalent.

	A parallel phenomenon arises in the present setting: filters (dually, ideals) on $\w$ which are inequivalent in the classical Kat\v{e}tov order may become equivalent in the Gamified Kat\v{e}tov order (e.g.\ MAD families). While category theory is certainly not needed to show this, Theorem~\ref{thm:mainthm} places it within a broader structural context. The Gamified Kat\v{e}tov order originates from the $\clt$-order, which is in order-reversing bijection with the lattice of subtoposes of $\Eff$ \cite[Theorem A.4.4.8]{Elephant}. In this sense, the identifications made by the $\glt$-order reflect, albeit indirectly, the category-theoretic equivalence used when defining the subtopos lattice of $\Eff$. Clarifying the relationship between these two notions of coarseness remains an interesting challenge, which we leave for future work.

	\bibliography{Eff}

\end{document}